\documentclass[letterpaper,12pt]{article}
\usepackage[margin=1in]{geometry}  
\usepackage{mathpazo}  
\usepackage{amsmath}
\usepackage{enumitem}
\usepackage{url}  
\usepackage{parskip}

\def\lf{\left\lfloor}
\def\rf{\right\rfloor}

\title{Finding the Year's Share in Day-of-Week Calculations}
\author{S. Kamal Abdali\thanks{
Email: {\em k.abdali@acm.org}}
}
\date{} 

\begin{document}
\maketitle

\begin{abstract}
The dominant part in the mental calculation of the day of the week for any given date is to determine 
the {\em year share}, that is, the contribution of the two-digit year part of the date. 
This paper describes a number of year share computation methods, some well-known and some new.  
The ``Parity Minus 3'' method, in particular, is a new alternative to the popular ``Odd+11'' method.
The paper categorizes the methods of year share computation, and presents simpler proofs of their 
correctness than usually provided. 
\\[6pt] \noindent \textbf{Keywords and phrases:} day of the week, calendar algorithms,  
doomsday method, first Sunday algorithm, mental arithmetic, year share 
\end{abstract}

\section{Introduction}
Finding the day of the week (DOW) for any given date is by now a trivial computational problem. 
While writing a program from scratch for computing the DOW is not difficult, most programming languages 
include libraries that provide routines for this purpose. 
Moreover, all existing computer operating systems and ``office'' applications have built-in facilities 
for this computation. 
For a discussion of several DOW algorithms best suited for programming, see the 
Wikipedia article ``Determination of the day of the week''\cite{wiki1}. 
A very comprehensive reference on DOW calculations in general is the German book 
{\em Enzyklop\"adie der Wochentagsberechnung} by Hans-Christian Solka\cite{solka}. 

While there is scarcely any need anymore for a new DOW algorithm for computers, the interest 
in {\em mental calculation} algorithms continues unabatedly. 
The best known, and undoubtedly one of the best, of such methods is the ``Doomsday Rule''\cite{berlekamp} 
invented in 1973 by John Conway. 
Another method with much merit, and for most people the easiest to use, is the 
``First Sunday Algorithm''\cite{goddard} introduced by Robert Goddard in 2009.\footnote{
Hans-Christian Solka has kindly informed me that methods equivalent to First Sunday were published 
by the German magician C.~Willmann in 1896, by E.~Rogent and W.W.~Durbin in 1927, 
and by some other authors later as part of {\em Dominical letter} research. 
The references, to which I don't have access, are in \cite{solka}. 
}
Such methods proceed by isolating and finding the contributions of the century, two-digit year, month, 
and day parts of the date, and then adding or subtracting these contributions to form a total. 
The total is an integer which represents the serial of the DOW being sought in some ordering of the days. 
(A common convention is to assign the numbers 0, 1, 2, \ldots, 6 to Sunday, Monday, Tuesday, 
\ldots, Saturday, respectively.) 
If the total happens to not be in the range 0, 1, 2, \ldots, 6, it is simply reduced modulo 7 to yield a value 
in that desired range. 
Note that since the total is formed by adding or subtracting various terms, each of those terms can also be 
individually reduced modulo 7 as soon as it is computed or at any stage during summation. 

The century and month part contributions typically amount to constants that one memorizes, and the 
contribution of the day part typically consists of counting weekdays forward or backward from a date 
determined by previous steps. 
It is the contribution of the two-digit year part that takes most time (the largest number of seconds!) 
in mental calculation. 
This is the case in both the Doomsday Rule and the First Sunday Algorithm. 

The DOW of any fixed date within a year (i.e., a month-day combination) is advanced one day 
by every common year and two days by every leap year. 
For a (two-digit) year \(y\) within a century, the DOW advances \(y\) days because of common years 
and \(\lf\frac{y}{4}\rf\) days because of leap years, making the total DOW advance from the century year 
to the year \(y\) equal to \(y+\lf\frac{y}{4}\rf \), that is, 
\begin{equation} 
\lf \frac{5y}{4}\rf   \label{eq-ys} 
\end{equation} 

We use the term {\em year share} to refer to the expression \eqref{eq-ys} that represents the contribution 
of the two-digit year part of a date to its DOW. 
Of course, any expression congruent modulo 7 to \eqref{eq-ys} serves our purpose equally well. 
In some methods the year share is subtracted while forming the total sum, so they compute a value which 
equals (or is congruent modulo 7 to) the negative of the expression \eqref{eq-ys}.   

When a new method of computing the year share is proposed, it has to be proven correct by showing that 
its result is congruent modulo 7 to the expression \eqref{eq-ys}, 
or to minus this expression when that's what the method claims to compute. 
While the correctness of many methods is quite obvious and hardly requires a proof, 
for some methods it is not immediate. 
But in the published account of these methods, the proofs are often missing or are 
unnecessarily complicated. 
We have tried to provide correctness proofs by very similar, simple arguments. 
It is hoped that by treating the methods in a uniform way the present approach 
will provide a better understanding of the methods and will be helpful in devising new methods 
and better variants of old methods. 

\section{Methods of Year Share Computation} 
A number of alternative methods of year share computations have appeared in the literature. 
These methods work essentially by substituting some expression for the variable \(y\) in \eqref{eq-ys} 
and transforming the resulting expression into one that is easy to evaluate by mental arithmetic. 
The transformations try to minimize the steps in the integer divisions by 4 and reduction modulo 7, 
and, of course, to maximize the work with small integers. 

We describe a selection of known methods here and also suggest a few new ones 
(Methods 2, 5c, 5f, 7, and 10). 
We have placed the methods into three categories of  (1) ``special'' methods that avoid division by 
a divisor other than 2, (2) division by integers larger than 2, and (3) operations on the individual digits 
of the two-digit year. 
The category (3) is really a subcategory of (2) because its methods involve division by 10. 
However, the methods in (3) require somewhat different mental arithmetic. 

\subsection{Special Methods} 
The overriding advantage of these methods is that you need to keep manipulating only a single variable 
(testing it, increasing or decreasing it, halving it). 
By contrast, the methods of the next sections require you to remember and work with  
several numbers. 

\begin{enumerate}[leftmargin=*] 

\item {\bf Odd+11 Method} 

This method, evolved from an idea proposed by Michael Walters in 2008, is described more formally 
in a 2011 paper by Chamberlain Fong and Michael Walters\cite{fong2}. 
Of all the present alternatives for year share computation, this method seems to be the quickest. 

The ``Odd+11'' operation takes a given year \(y\) (between 0 and 99, inclusive) as input 
and produces the negative of the year share as output, by proceeding as follows: 
\begin{enumerate}[label=\roman*.] 
\item 
Set the value of \(YS\) to that of \(y\). In symbols, \(YS \leftarrow y\). 
\item 
If \(YS\) is odd, then increase it by 11, i.e., \(YS \leftarrow YS+11\), else leave it unchanged. 
\item 
Halve \(YS\), i.e., \(YS \leftarrow YS/2\). 
\item 
If \(YS\) is odd, increase it by 11, i.e., \(YS \leftarrow YS+11\), else leave it unchanged. 
\end{enumerate} 
Fong and Walter\cite{fong2} add two more steps of doing \(YS \leftarrow YS \bmod 7\) and 
then \(YS \leftarrow 7-YS\) to turn the result into a positive year share. 
We omit these steps since the negative year share is what the ``First Sunday Algorithm'' needs 
anyway, and, moreover, any result in the mod 7 congruence class is acceptable. 

To show that the result of applying ``Odd+11'' to the input \(y\) is congruent modulo 7 
to the negative of \eqref{eq-ys}, we first express \(y\) as a polynomial as follows: 
Divide \(y\) by 4; call the quotient \(a\); 
divide the remainder (whose value is between 0 and 3, inclusive) by 2; 
call the new quotient \(b\) and the new remainder \(c\). 
We can now write 
\begin{equation} 
y = 4a + 2b +c, \label{eq-od1} 
\end{equation} 
where \(0 \leq b \leq 1 \ \textrm{and}\ 0 \leq c \leq 1.\) 

Let us first apply the steps of ``Odd+11'' to this value of \(y\). 
\begin{enumerate}[label=\roman*.] 
\item 
\[YS = 4a+2b+c.\] 
\item 
If \(YS\) is odd, i.e., if \(c=1\), then increase \(YS\) by 11, else leave \(YS\) unchanged. 
This is the same as writing 
\[YS=4a+2b+12c.\] 
\item 
Halve \(YS\), i.e., 
\[YS=2a+b+6c.\] 
\item 
If \(YS\) is odd, i.e., if \(b=1\), then increase \(YS\) by 11, else leave \(YS\) unchanged. 
This is the same as writing 
\begin{equation} 
YS=2a+12b+6c.  \label{eq-od2} 
\end{equation} 
\end{enumerate} 
Next, let us evaluate the negative of the expression \eqref{eq-ys} for the value of \(y\) 
given by \eqref{eq-od1}. 
\begin{equation} 
 - \lf\frac{5y}{4}\rf    
=  - \lf\frac{20a+10b+5c}{4}\rf   
=  -\left(5a + 2b +c +\lf\frac{2b +c}{4}\rf \right) 
=  -5a -2b -c   \label{eq-od3} 
\end{equation} 
As \eqref{eq-od2} and \eqref{eq-od3} differ by a multiple of 7 (viz. \(7a+14b+ 7c\)), 
they are congruent modulo 7. 

Fong and Walters's own proof in \cite{fong2} is longer and more complicated. 

\item{\bf Parity Minus 3 Method} 

This new method is inspired by and very similar to ``Odd+11'', but involves  
smaller integers. 
The ``Parity Minus 3'' operation takes a given year \(y\) (between 0 and 99, inclusive) as input 
and produces the negative of the year share as output, by proceeding as follows: 
\begin{enumerate}[label=\roman*.] 
\item 
Set the value of \(YS\) to that of \(y\). In symbols, \(YS \leftarrow y\).  
\item 
Check and remember \(YS\)'s parity (odd or even). 
If \(YS\) is odd, then decrease it by 3, i.e., \(YS \leftarrow YS-3\), else leave it unchanged. 
\item 
Halve \(YS\), i.e., \(YS \leftarrow YS/2\). 
\item 
If \(YS\)'s parity (odd or even) has changed, decrease \(YS\) by 3, i.e., \(YS \leftarrow YS-3\), 
else leave it unchanged. 
\end{enumerate} 

Examples: 
\begin{enumerate}[label=(\alph*)] 
\item \(y=24\) 

Step 1: \(YS=24\). Step 2: Even, hence \(YS=24\). Step 3: Halve, so \(YS=12\). 
Step 4: Even, so parity unchanged, hence answer is \(YS=12\). 
\item \(y=37\) 

Step 1: \(YS=37\). Step 2: Odd, hence \(YS=37-3=34\). Step 3: Halve, so \(YS=17\). 
Step 4: Odd, so parity unchanged, hence answer is \(YS=17\). 
\item \(y=58\) 

Step 1: \(YS=58\). Step 2: Even, hence \(YS=58\). Step 3: Halve, so \(YS=29\). 
Step 4: Odd, so parity changed, hence answer is \(YS=29-3=26\). 
\item \(y=79\) 

Step 1: \(YS=79\). Step 2: Odd, hence \(YS=79-3=76\). Step 3: Halve, so \(YS=38\). 
Step 4: Even, so parity changed, hence answer is \(YS=38-3=35\). 
\end{enumerate} 

To show that the result of applying ``Parity Minus 3'' to the input \(y\) is congruent modulo 7 
to the negative of \eqref{eq-ys}, we first express \(y\) as a polynomial as follows: 
Divide \(y\) by 4; call the quotient \(a\); 
divide the remainder (whose value is between 0 and 3, inclusive) by 2; 
call the new quotient \(b\) and the new remainder \(c\). 
We can now write 
\begin{equation} 
y = 4a + 2b +c, \label{eq-par1} 
\end{equation} 
where \(0 \leq b \leq 1 \ \textrm{and}\ 0 \leq c \leq 1.\) 

Let us first apply the steps of ``Parity Minus 3'' to this value of \(y\). 
\begin{enumerate}[label=\roman*.] 
\item 
\[YS = 4a+2b+c.\] 
\item 
If \(YS\) is odd, i.e., if \(c=1\), then decrease \(YS\) by 3, else leave \(YS\) unchanged. 
This is the same as writing 
\[YS=4a+2b-2c.\] 
\item 
Halve \(YS\), i.e., 
\[YS=2a+b-c.\] 
\item 
Since each of \(b\) and \(c\) is either 0 or 1, the new parity of \(YS=2a+b-c\) is odd or even depending, 
respectively, on whether \(b\ne c\) or \(b=c\). 
The old parity of \(YS\) determined in Step 2 of the method was odd or even depending on whether 
\(c=1\) or \(c=0\). 
Thus the new parity is different from or same as the old one depending on whether \(b=1\) or \(b=0\). 
\(YS\) is to be decreased by 3 if the parity has changed, i.e., \(b=1\), and is to be left unchanged if \(b=0\). 
This is the same as writing 
\begin{equation} 
YS=2a-2b-c.  \label{eq-par2} 
\end{equation} 
\end{enumerate} 
Next, let us evaluate the negative of the expression \eqref{eq-ys} for the value of \(y\)  
given by \eqref{eq-par1}. 
\begin{equation} 
- \lf\frac{5y}{4}\rf   
=   - \lf\frac{20a+10b+5c}{4}\rf    
=   -\left(5a + 2b +c +\lf\frac{2b +c}{4}\rf \right) 
=  -5a -2b -c   \label{eq-par3} 
\end{equation} 
As \eqref{eq-par2} and \eqref{eq-par3} differ by a multiple of 7, they are congruent modulo 7. 

\end{enumerate} 

\subsection{Division by Various Integers} 

\begin{enumerate}[resume,leftmargin=*] 
\item{\bf Division by 12} 

In the original version of the Doomsday Rule\cite{berlekamp}, Conway states the method 
(originally due to Lewis Carroll, see Gardner\cite{gardner}) to compute the year 
share as follows: 

``add the number of {\em dozens} [\ldots\  in \(y\)], the {\em remainder} after  
[\ldots\  the dozens are taken out], and the number of {\em fours} in the {\em remainder}''. 

That is, the year share is 
\begin{equation}  
\lf\frac{y}{12}\rf + y \bmod 12 + \lf\frac{y \bmod 12}{4}\rf. \label{eq-c1} 
\end{equation} 

To prove that the method computes the year share correctly, let us write 
\[y=12q+r, \ \textrm{where}\  0 \leq r < 12.\] 
Now \eqref{eq-c1} can be rewritten as 
\begin{equation} 
q+r+\lf\frac{r}{4}\rf    \label{eq-c2} 
\end{equation} 
With \(12q+r\) substituted for \(y\) in \eqref{eq-ys}, year share equals 
\begin{equation} 
\lf\frac{5y}{4}\rf = \lf\frac{60q+5r}{4}\rf = 15q+r+\lf\frac{r}{4}\rf.  \label{eq-c3} 
\end{equation} 
Since \eqref{eq-c2} and \eqref{eq-c3} differ by a multiple of 7, they are congruent modulo 7.

\item {\bf Division by 4} 

This method called {\em Highest Multiple of Four} by YingKing Yu\cite{yu} uses a very simple, easy to 
remember calculation to produce the negative of the year share. 
Solka\cite{solka} credits Carl Willmann with a much earlier equivalent method. 
One thinks of the two-digit year as the sum of a multiple of 4 and a remainder which 
is 0, 1, 2, or 3. Then the method computes ``half of that multiple of 4, minus the remainder''. 
Incidentally, for some people ``the closet multiple of 4 not larger than the given year'' is easier to remember 
as the current or previous ``leap year'', ``Olympics year'', or the ``US Presidential Election year''. 

To prove that the method works correctly, let us write \(y=4q+r, \ \textrm{where}\  0 \leq r < 4\). 
Then the output of the method is \(2q-r\). 
That the value of the negative of the expression \eqref{eq-ys} when \(4q+r\) is substituted  
for \(y\) is congruent modulo 7 to the expression \(2q-r\)  is shown as follows: 
\begin{equation*}  
- \lf\frac{5y}{4}\rf = - \lf\frac{20q+5r}{4}\rf = -(5q+r+\lf\frac{r}{4}\rf) = -5q-r \equiv 2q-r \pmod 7 
\end{equation*}

\item{\bf Division by Other Integers} 

Divisions by 4 and 12 furnish nice methods, as we have seen above, 
because these divisions result in formulas that require simple arithmetic. 
Division by 10 also has nice properties and will be covered in the next section. 
The divisors 5, 11, 16, and 17 also lead to simple formulas for the year share. 
Solka\cite{solka} gives formulas for several divisors, including 8, 12, 16, 20, and 24, 
grouped together into a ``universal approach'' section. 

In general, for a divisor \(d\), we write \(y=dq+r\), where \(0 \le r < d\), and evaluate 
the year share expression \eqref{eq-ys}, i.e., \(\lf\frac{5dq+5r}{4}\rf\) 
(or its negative). 
After reducing \(5d\) modulo 28, we can expand the expression into the form 
\(aq+r + \lf\frac{bq+r}{4}\rf\) for some integers \(a\) and \(b\).  
We can further play with this expression in various ways, e.g., increase \(a\) by any 
integer \(k\) and compensate for that change by decreasing \(b\) by \(4k\). 
Below we show only those divisors \(d\) between 5 and 20 for which \(a\) and \(b\) have values 
\(0, +1\), or \(-1\), because any other value would require extra multiplications that would 
complicate mental arithmetic.  
For the sake of comparison, we include the formulas for divisors 4 and 12 
given in the previous sections. 

\begin{enumerate}[label=(\alph*)] 
\item \(d=4\): 
\begin{equation*} 
\textrm{Negative year share} = 2q-r. 
\end{equation*} 
    
\item \(d=5\): 
\begin{equation*} 
\textrm{Negative year share} = q-r-\lf \frac{q+r}{4}\rf. 
\end{equation*} 
    
\item \(d=11\): 
\begin{equation*} 
\textrm{Positive year share} = r+\lf \frac{r-q}{4}\rf. 
\end{equation*} 
Note that the dividend \(r-q\) in the expression \(\lf\frac{r-q}{4}\rf\) can be negative. 
The evaluation of integer quotients in such cases requires some extra care.  
This is discussed at the beginning of Section~\ref{sec:digits}. 
 
\item \(d=12\): 
\begin{equation*} 
\textrm{Negative year share} = q+r+\lf \frac{r}{4}\rf. 
\end{equation*} 
    
\item \(d=16\): 
\begin{equation*} 
\textrm{Positive year share} = -q+r+\lf \frac{r}{4}\rf. 
\end{equation*} 
    
 \item \(d=17\): 
\begin{equation*} 
\textrm{Positive year share} = r+\lf \frac{q+r}{4}\rf. 
\end{equation*} 
\end{enumerate} 
  
\end{enumerate} 

\subsection{Methods Operating on the Year's Individual Digits} 
\label{sec:digits}  
The advantage of these methods is that they involve arithmetic with smaller numbers than those arising 
in the methods of the previous category, 
and dividing these numbers by 4 or reducing them modulo 7 is quite easy to do mentally. 

The earliest such method is credited by Solka\cite{solka} to L.T.~Sakharovski. 
This method, published in 1957, has assigned codes to the tens and 
units digits of the year, and these codes are added together to get the year share. 
As our interest is mainly in year share {\em computation}, we will not reproduce 
Sakharovski's table (given in \cite{solka}), 
and will describe five methods that do this computation in various ways. 

Method \ref{item:E} is different from the other methods in this section because it operates on the 
tens and units digits not of the two-digit year but of the highest multiple of four not exceeding the year. 
Method \ref{item:Aa} seems to be the simplest, but in Method \ref{item:Ab} 
we work with only one number at a time (while keeping one sign in memory), similarly to ``Odd+11'' 
or ``Parity Minus 3''. 

Some of the methods in the present category involve integer divisions with the dividends 
allowed to be negative. 
(However, the divisor is required to be positive.) 
Our needed congruences will hold only if we define the integer quotient in such a division as follows: 
\begin{equation} 
\lf\frac{p}{q}\rf = 
\begin{cases} 
     \lf\frac{p}{q}\rf,  & \textrm{if }  p\geq 0 \ (\textrm{and } q > 0) \nonumber \\ 
     \ \nonumber \\ 
   - \lf\frac{|p|}{q}\rf, & \textrm{if } p < 0 \textrm{ and } |p| \bmod q = 0 \ (\textrm{and } q > 0) \nonumber \\ 
     \ \nonumber \\ 
   - \left( \lf\frac{|p|}{q}\rf + 1\right), & \textrm{if }  p < 0 \textrm{ and } |p| \bmod q > 0 \ (\textrm{and } q > 0)  
\end{cases} 
\end{equation} 

Methods \ref{item:Aa} through \ref{item:Ab} have much in common and are essentially variations on the same theme.  
Suppose the tens and units digits of the two digit-year part are, respectively, \(t\) and \(u\). 
That is \(y=10t+u\). 
Let's evaluate the negative of the year share expression \eqref{eq-ys} in terms of \(t\) and \(u\): 
\begin{equation} 
 -  \lf\frac{5y}{4}\rf\  
= - \left(\lf\frac{50t+5u}{4}\rf \right)  \label{eq-TC1}  
\end{equation} 
Any expression derived by adding a multiple of 28 to the numerator of the fraction in \eqref{eq-TC1} is 
obviously congruent modulo 7 to \eqref{eq-TC1}, and is hence just another expression for the negative 
year share. 
Examples are: 
\begin{align} 
& - \left(\lf\frac{22t + 5u}{4} \rf \right)  \label{eq-TC2} \\ 
& - \left(\lf\frac{-6t + 5u}{4} \rf \right) \label{eq-TC3} 
\end{align} 
Method  \ref{item:Aa}, \ref{item:F}, and \ref{item:W} are derived from \eqref{eq-TC2} with or without 
the negative sign, and Method \ref{item:Ab} in essence evaluates \eqref{eq-TC3} directly.

\begin{enumerate}[resume,leftmargin=*] 
\item {\bf Computing with Year's Individual Digits (Eisele)}\label{item:E} 

The following method by Martin Eisele (citation in \cite{solka}) is unique as it operates on the tens 
and units digits not of the given year but of the largest multiple of four not exceeding that year. 
Let \(y\) be the two-digit year part of the date.  
Let \(q\) and \(r\) be the integer quotient and remainder when \(y\) is divided by 4, that is, 
\begin{equation} 
y=4q+r, \ \textrm{where}\  0 \le r < 4.   \label{eq-tuE1} 
\end{equation} 
Let \(t\) and \(u\) be the tens and units digits of \(4q\), that is, 
\begin{equation} 
4q=10t+u, \ \textrm{where}\  0 \le u < 10. \label{eq-tuE2} 
\end{equation} 
Then, according to Eisele, the year share is  
\begin{equation} 
2t-\frac{u}{2}+r.    \label{eq-tuE3} 
\end{equation} 
Note that \(u\) must be even since by  \eqref{eq-tuE2}  \(10t+u\) is a multiple of 4. 
Thus we have 
\begin{equation} 
2q=5t+\frac{u}{2}, \ \textrm{where}\  u = 0, 2, 4, 6, \textrm{or}\  8.   \label{eq-tuE4} 
\end{equation} 
To show the correctness of Eisele's method, we verify that the year share expression \eqref{eq-ys} 
with \(4q+r\) substituted for \(y\) is congruent modulo 7 to \eqref{eq-tuE3}: 
\begin{align} 
& \lf\frac{5y}{4}\rf    
= \lf\frac{20q+5r}{4}\rf, \ \textrm{by}\  \eqref{eq-tuE1}  \nonumber \\ 
= & \:    5q+r + \lf\frac{r}{4}\rf   = 5q+r  
\equiv  -2q+r \pmod 7 
=  -5t-\frac{u}{2}+r, \ \textrm{by}\  \eqref{eq-tuE4} \nonumber \\ 
\equiv & \:  2t-\frac{u}{2}+r \pmod 7   \nonumber 
\end{align} 

A method by Alexander Harringer (see \cite{solka}) turns out to be a variation 
of the above method. 
For computing the year share, Harringer  proposes the formula 
\begin{equation} 
2t +3u +r.   \label{eq-tuE5} 
\end{equation} 
instead of Eisele's \eqref{eq-tuE3}. 
Notice that \eqref{eq-tuE5} and \eqref{eq-tuE3} differ by the quantity \(\frac{7u}{2}\). 
As \(u\) is an even number, this quantity is a multiple of 7, and hence the formulas by 
Harringer and Eisele are congruent modulo 7.

\item {\bf Computing with Year's Individual Digits (Aa)}\label{item:Aa} 

Let the tens digit and units digit in the year part be \(t\) and \(u\), respectively. 
That is, \(y=10t+u\). 
Then the negative of the year share can be found in this way: 
\begin{enumerate}[label=\roman*.] 
\item 
Compute \( \lf \frac{2t+u}{4}\rf \).  
\item 
Add \(u\) to above.  
\item 
Subtract the above sum from \(2t\).\\ 
This is the year share. 
\end{enumerate} 
Example. Suppose  \(year\) is 59. Then: 
\begin{enumerate}[label=\roman*.] 
\item 
From \(t=5\) and \(u=9\), we compute \(2t+u=19\) and then its quarter which is 4. 
\item 
Adding \(u=9\) to it, we get 13. 
\item 
Subtracting this from \(2t=10\), we get \(10-13 = -3\).  
This is the negative year share. 
We can reduce it modulo 7 immediately to 4 or just leave it as \(-3\) 
to be reduced modulo 7 in a later step of DOW calculation. 
\end{enumerate} 
The result of performing the steps of the method with input digits \(t\) and \(u\) is 
\begin{equation} 
2t - \left(\lf \frac{2t+u}{4}\rf + u \right) \label{eq-tuAa1} 
\end{equation} 
To show that this is the negative year share, we evaluate the negative of expression \eqref{eq-ys} 
with \(10t+u\) substituted for \(y\). 
\begin{align} 
& - \lf\frac{5y}{4}\rf    
= - \lf\frac{50t+5u}{4}\rf  
= - \left(12t+u + 2t + \lf\frac{2t+u}{4}\rf \right)  \nonumber \\ 
= &  -12t - \left(\lf\frac{2t+u}{4}\rf + u \right)  \label{eq-tuAa2} 
\end{align} 
Since \eqref{eq-tuAa1} and \eqref{eq-tuAa2} differ by a multiple of 7, they are congruent modulo 7.

\item {\bf Computing with Year's Individual Digits (Fong)} \label{item:F} 

The following method by Chamberlain Fong\cite{fong1} operates on the year's individual digits and computes 
the positive year share. 
(Fong also credits YingKing Yu with this method, and cites Yu's work in \cite{fong1}.) 
Let \(t\) and \(u\) be, respectively, the tens and units digits of the year. 
Then the method computes the year share as 
\begin{equation} 
2t + 10 (t \bmod 2) +u + \lf\frac{2 (t \bmod 2) + u}{4}\rf.  \label{eq-tuF1} 
\end{equation} 
To prove that the year share so computed is correct, we proceed as follows. 
Substituting \(t=2t_1+t_2\), where \(0 \le t_2 \le 1\), we write \eqref{eq-tuF1} as 
\begin{equation} 
 2(2 t_1 + t_2) + 10 t_2 +u + \lf\frac{2 t_2 + u}{4}\rf   
 =   4 t_1 + 12 t_2 + u + \lf\frac{2 t_2 + u}{4}\rf.  \label{eq-tuF2} 
\end{equation} 
Substituting \(t=2 t_1 + t_2\), hence \(y=10t+u=20t_1 + 10t_2+u\), we evaluate the expression \eqref{eq-ys} 
for year share as follows: 
\begin{equation} 
  \lf\frac{5y}{4}\rf   
=   \lf \frac{100 t_1 + 50 t_2 + 5u}{4}\rf   
 =  25t_1 + 12 t_2 + u + \lf \frac{2 t_2 + u}{4}\rf   \label{eq-tuF3} 
\end{equation} 
Since \eqref{eq-tuF2} and \eqref{eq-tuF3} differ by a multiple of 7, they are congruent modulo 7. 

\item {\bf Computing with Year's Individual Digits (Wang)} \label{item:W} 

This method by Xiang-Sheng Wang\cite{wang} computes the positive year share. 
Let \(t\) and \(u\) be, respectively, the tens and units digits of the year. 
Then the method computes the year share as 
\begin{equation} 
u-t  + \lf \frac{u}{4} - \frac{t}{2} \rf.  \label{eq-tuW1} 
\end{equation} 
To prove that the year share so computed is correct,  we evaluate the negative of the year share 
expression \eqref{eq-ys} with \(10t+u\) substituted for \(y\). 
\begin{align} 
& \lf\frac{5y}{4}\rf    
=    \lf\frac{50t + 5u}{4}\rf  
=  \lf \frac{4u+52t+u-2t}{4}\rf  
=  u+13t + \lf\frac{u-2t}{4}\rf  \nonumber \\ 
= & \: u+13t + \lf \frac{u}{4} - \frac{t}{2}\rf.  \label{eq-tuW2} 
\end{align}
Since \eqref{eq-tuW1} and \eqref{eq-tuW2} differ by a multiple of 7, they are congruent modulo 7.  

Note that the fraction in \eqref{eq-tuW1} can be negative, so its floor has to be evaluated carefully  
following  the procedure stated at the beginning of Section~\ref{sec:digits}.

\item {\bf Computing with Year's Individual Digits (Ab)}\label{item:Ab} 

Here is another method that operates on the individual digits of a two-digit year. 
Let the tens digit and units digit in the year part be \(t\) and \(u\), respectively; 
that is, \(y=10t+u\). 
Then the negative year share is found in this way: 
\begin{enumerate}[label=\roman*.] 
\item 
Compute \(5u-6t\), and let its absolute (positive) value be  \(a\). 
Also remember its sign (`plus' or `minus'). 
\item 
Compute \(b =  \lf\frac{a}{4}\rf \). 
If there was a non-zero remainder, and the sign in the previous step was `minus', increase \(b\) by 1. 
\item 
Affix the opposite sign to \(b\). 
(That is, make it \(-b\) if the sign was `plus', and \(+b\) i.e., just \(b\) if the sign was `minus'. ) 
This result is the negative year share. \\ 
If this value turns out to be negative or larger than 6, we can reduce it modulo 7 immediately or just leave it as is 
to be reduced modulo 7 in a later step of DOW calculation. 
\end{enumerate} 
Example. Suppose  \(year\) is 87. Then: 
\begin{enumerate}[label=\roman*.] 
\item 
From \(t=8\) and \(u=7\), we compute 
\(5u-6t=5\times 7 - 6\times 8 = 35-48 = -13\). 
So we have \(a=13\) and sign \(=\) `minus'. 
\item 
By dividing \(a\) by 4, we get \(b=3\). Since there was a nonzero remainder in the division, 
and the remembered sign is `minus', we add 1 to \(b\), making \(b=4\). 
\item 
Since the remembered sign is `minus', we attach the opposite sign to \(b\), making it \(+4\), i.e., 4. 
So the negative year share is 4. 
\end{enumerate} 
The result of performing the steps of the method with input digits \(t\) and \(u\) is 
\begin{equation} 
- \lf\frac{5u-6t}{4}\rf  \label{eq-tuAb1} 
\end{equation} 
To show that this is the negative year share, we evaluate the negative of expression \eqref{eq-ys} 
with \(10t+u\) substituted for \(y\). 
\begin{equation} 
 - \lf\frac{5y}{4}\rf    
= - \lf\frac{50t+5u}{4}\rf   
= - \lf\frac{5u-6t+56t}{4}\rf   
\equiv  - \left(\lf\frac{5u - 6t}{4} \rf  + 14t \right)    \label{eq-tuAb2} 
\end{equation} 
Since \eqref{eq-tuAb1} and \eqref{eq-tuAb2} differ by a multiple of 7, they are congruent modulo 7.

\end{enumerate} 

\section{Concluding Remarks} 
This article is concerned only with the year share part of DOW calculation, not with any other details of 
the DOW computation methods. 
The year share part is where most of the calculation time is spent. 
We have tried to describe the methods with a uniform, systematic approach, and have 
provided simple proofs of their correctness. 

\section{Acknowledgment} 
I am thankful to Hans-Christian Solka, Robert Goddard, and Chamberlain Fong for reading an earlier 
draft of this article very carefully and making valuable suggestions.

\begin{flushleft} 
 
\end{flushleft}

\end{document}